\documentclass[a4paper,11pt]{amsart}
\usepackage{amsfonts,amsmath,amssymb,amsthm,graphicx}
\usepackage[lmargin=34mm,rmargin=34mm,tmargin=28mm,bmargin=28mm,footskip=10mm]{geometry}
\usepackage[sort&compress,numbers]{natbib}

\newtheorem{theorem}{Theorem}
\newtheorem{corollary}{Corollary}
\newtheorem{lemma}{Lemma}

\renewcommand{\baselinestretch}{1.1}
\newcommand{\corref}[1]{Corollary~\ref{cor:#1}}
\newcommand{\corlabel}[1]{\label{cor:#1}}
\newcommand{\lemlabel}[1]{\label{lem:#1}}
\newcommand{\lemref}[1]{Lemma~\ref{lem:#1}}
\newcommand{\thmref}[1]{Theorem~\ref{thm:#1}}
\newcommand{\thmlabel}[1]{\label{thm:#1}}
\newcommand{\figref}[1]{Figure~\ref{fig:#1}}
\newcommand{\figlabel}[1]{\label{fig:#1}}

\newcommand{\twolemref}[2]{Lemmata~\ref{lem:#1} and \ref{lem:#2}}

\newcommand{\Figure}[4][htb]{
\begin{figure}[#1]
	\begin{center}#3\end{center}
	\caption{\figlabel{#2}#4}
\end{figure}
}

\newcommand{\Oh}[1]{\ensuremath{\protect\mathcal{O}(#1)}}
\newcommand{\bracket}[1]{\ensuremath{\protect\left(#1\right)}}

\newcommand{\half}{\ensuremath{\protect\tfrac{1}{2}}}
\newcommand{\bfrac}[2]{\ensuremath{\protect\bracket{\frac{#1}{#2}}}}

\newcommand{\etal}{~et~al.~}

\newcommand{\e}{\textsf{\textup{e}}}

\begin{document}

\date{\today}
\title[Graphs with Large Queue-Number]{Bounded-Degree Graphs have\\ Arbitrarily Large Queue-Number}
\author{David R.~Wood}
\address{Departament de Matem{\'a}tica Aplicada II, Universitat Polit{\`e}cnica de Catalunya, Barcelona, Spain}
\email{david.wood@upc.edu}
\thanks{Supported by the Government of Spain grant MEC SB2003-0270, and by the projects MCYT-FEDER BFM2003-00368 and Gen.\ Cat 2001SGR00224.}


\keywords{graph, queue layout, queue-number}

\subjclass[2000]{05C62 (graph representations)}

\begin{abstract} 
It is proved that there exist graphs of bounded degree with arbitrarily large queue-number. In particular, for all $\Delta\geq3$ and for all sufficiently large $n$, there is a simple $\Delta$-regular $n$-vertex graph with queue-number at least $c\sqrt{\Delta}n^{1/2-1/\Delta}$ for some absolute constant $c$.
\end{abstract}


\maketitle

\section{Introduction}

We consider graphs possibly with loops but with no parallel edges. A graph without loops is \emph{simple}. Let $G$ be a graph with vertex set $V(G)$ and edge set $E(G)$. If $S\subseteq E(G)$ then $G[S]$ denotes the spanning subgraph of $G$ with edge set $S$. We say $G$ is \emph{ordered} if $V(G)=\{1,2,\dots,|V(G)|\}$. Let $G$ be an ordered graph. Let $\ell(e)$ and $r(e)$ denote the endpoints of each edge $e\in E(G)$ such that $\ell(e)\leq r(e)$. Two edges $e$ and $f$ are \emph{nested} and $f$ is \emph{nested inside} $e$ if $\ell(e)<\ell(f)$ and $r(f)<r(e)$. An ordered graph is a \emph{queue} if no two edges are nested. Observe that the left and right endpoints of the edges in a queue are in first-in-first-out order---hence the name `queue'. An ordered graph $G$ is a $k$-\emph{queue} if there is a partition $\{E_1,E_2,\dots,E_k\}$ of $E(G)$ such that each $G[E_i]$ is a queue. 

Let $G$ be an (unordered) graph. A \emph{$k$-queue layout} of $G$ is a $k$-queue that is isomorphic to $G$. The \emph{queue-number} of $G$ is the minimum integer $k$ such that $G$ has a $k$-queue layout. Queue layouts and queue-number were introduced by Heath\etal\citep{HLR-SJDM92,HR-SJC92} in 1992, and have applications in sorting permutations \citep{Tarjan72a, EI71, OS-CN93, IN-IPJ77, Pratt73}, parallel process scheduling \citep{BCLR-JPDC96}, matrix computations \citep{Pemmaraju-PhD}, and graph drawing \citep{DMW-SJC05, GLM-CGTA05}. Other aspects of queue layouts have been studied in \citep{Wood-Queue-DMTCS05, DujWoo-DMTCS04, RM-COCOON95, SS00, DujWoo-DMTCS05, DPW-DMTCS04, Ganley-HalinGraphs-95}. 

Prior to this work it was unknown whether graphs of bounded degree have bounded queue-number. The main contribution of this note is to prove that there exist graphs of bounded degree with arbitrarily large queue-number. 

\begin{theorem}
\thmlabel{QueueDegree}
For all $\Delta\geq3$ and for all sufficiently large $n>n(\Delta)$, there is a simple $\Delta$-regular $n$-vertex graph with queue-number at least $c\sqrt{\Delta}n^{1/2-1/\Delta}$ for some absolute constant $c$.
\end{theorem}

The best known upper bound on the queue-number of a graph with maximum degree $\Delta$ is $\e(\Delta n/2)^{1/2}$ due to \citet{DujWoo-DMTCS04} (where \e\ is the base of the natural logarithm). Observe that for large $\Delta$, the lower bound in \thmref{QueueDegree} tends toward this upper bound. Although for specific values of $\Delta$ a gap remains. For example, for $\Delta=3$ we have an existential lower bound of $\Omega(n^{1/6})$ and a universal upper bound of \Oh{n^{1/2}}.

\section{Proof of \thmref{QueueDegree}}

The proof of \thmref{QueueDegree} is modelled on a similar proof by \citet{BMW-EJC06}. Basically, we show that there are more graphs $\Delta$-regular graphs than graphs with bounded queue-number.  
The following lower bound on the number of $\Delta$-regular graphs is a corollary of more precise bounds due to \citet{BC-JCTA78}, \citet{Wormald78}, and \citet{McKay-AC85}; see \citep{BMW-EJC06}.

\begin{lemma}[\citep{BC-JCTA78,Wormald78,McKay-AC85}]
\lemlabel{NumberRegular}
For all integers $\Delta\geq1$ and for sufficiently large $n\geq n(\Delta)$, 
the number of labelled simple $\Delta$-regular $n$-vertex graphs is at least 
$$\bfrac{n}{3\Delta}^{\Delta n/2}$$
\end{lemma}

It remains to count the graphs with bounded queue-number. We will need the following two lemmas from the literature, whose proofs we include for completeness. A \emph{rainbow} in an ordered graph is a set of pairwise nested edges. 

\begin{lemma}[\citep{HR-SJC92,DujWoo-DMTCS04}]
\lemlabel{Rainbow}
An ordered graph $G$ is a $k$-queue if and only if $G$ has no $(k+1)$-edge rainbow. 
\end{lemma}

\begin{proof}
The necessity is obvious. For the sufficiency, suppose $G$ has no $(k+1)$-edge rainbow. For every edge $e$ of $G$, if $i-1$ edges are pairwise nested inside $e$, then assign $e$ to the $i$-th queue. 
\end{proof}

\begin{lemma}[\citep{DujWoo-DMTCS04}]
\label{lem:Edges}
Every ordered $n$-vertex graph with no two nested edges has at most $2n-1$ edges.
\end{lemma}

\begin{proof}
If $v+w=x+y$ for two distinct edges $vw$ and $xy$, then $vw$ and $xy$ are nested. The result follows since $2\leq v+w\leq 2n$. 
\end{proof}

Let $g(n)$ be the number of queues on $n$ vertices. To bound $g(n)$ we adapt a proof of a more general result by \citet{Klazar00}; also see \citep{Tardos-JCTA05, MT-JCTA04} for other related and more general results. 

\begin{lemma}
\lemlabel{NumberNested}
$g(n)\leq c^n$ for some absolute constant $c$.
\end{lemma}

\begin{proof}
Say $G$ is an ordered $n$-vertex graph. Let $G'$ be an ordered $2n$-vertex graph obtained by the following \emph{doubling} operation. For every edge $vw$ of $G$, add to $G'$ a nonempty set of edges between $\{2v-1,2v\}$ and $\{2w-1,2w\}$, no pair of which are nested. 

Every ordered $2n$-vertex ordered graph with no two nested edges can be obtained from some ordered $n$-vertex graph with no two nested pair of edges by a \emph{doubling} operation. To see this, merge every second pair of vertices, introduce a loop between merged vertices, and replace any resulting parallel edges by a single edge. The ordered graph that is obtained has no nested pair of edges. 

In the doubling operation, there are $11$ possible ways to add a nonempty set of non-nested edges between $\{2v-1,2v\}$ and $\{2w-1,2w\}$, as illustrated in \figref{eleven}. Thus $g(2n) \leq 11^{2n-1}\cdot g(n)$, since $G'$ has at most $2n-1$ edges by \lemref{Edges}. It follows that $g(n)\leq 11^{2n}$.
\end{proof}

\Figure{eleven}{\includegraphics{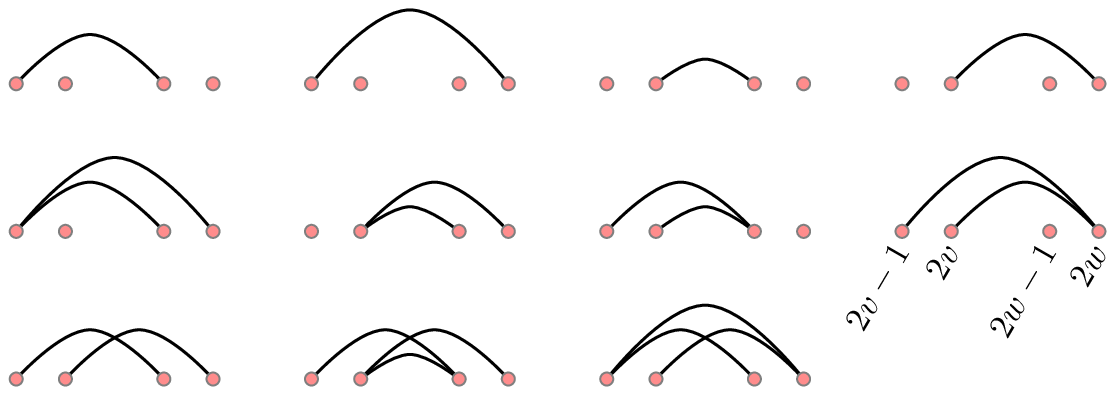}}{The $11$ possible ways to add a nonempty set of non-nested edges between $\{2v-1,2v\}$ and $\{2w-1,2w\}$.}


\twolemref{Rainbow}{NumberNested} imply the following.

\begin{corollary}
\corlabel{NumberRainbow}
The number of $k$-queues on $n$ vertices is at most $c^{kn}$ for some absolute constant $c$.
\end{corollary}

It is easily seen that \lemref{NumberRegular} and \corref{NumberRainbow} imply a lower bound of $c(\Delta/2-1)\log n$ on the queue-number of some $\Delta$-regular $n$-vertex graph. To improve this logarithmic bound to polynomial, we now give a more precise analysis of the number of $k$-queues that also accounts for the number of edges in the graph.


Let $g(n,m)$ be the number of $k$-queues on $n$ vertices and $m$ edges.

\begin{lemma}
\label{lem:NumberEdges}
$$g(n,m)\leq\begin{cases}
\binom{n}{2m}\cdot c^{2m}	&\textup{, if } m\leq\frac{n}{2}\\
c^{n}						&\textup{, if } m>\frac{n}{2},
\end{cases}$$
for some absolute constant $c$.
\end{lemma}

\begin{proof}
By \lemref{NumberNested}, we have the upper bound of $c^n$ regardless of $m$. Suppose that $m\leq\frac{n}{2}$. An $m$-edge graph has at most $2m$ vertices of nonzero degree. Thus every $n$-vertex $m$-edge queue is obtained by first choosing a set $S$ of $2m$ vertices, and then choosing a queue with $|S|$ vertices. The result follows.
\end{proof}

Let $g(n,m,k)$ be the number of $k$-queues on $n$ vertices and $m$ edges.

\begin{lemma}
\lemlabel{NumberRainbowEdges}
For all integers $k$ such that $\frac{2m}{n}\leq k\leq m$,  $$g(n,m,k)\leq\bfrac{ckn}{m}^{2m}$$ for some absolute constant $c$.
\end{lemma}

\begin{proof}
Fix nonnegative integers $m_1\leq m_2\leq\dots\leq m_k$ such that $\sum_im_i=m$. Let $g(n;m_1,m_2,\dots,m_k)$ be the number of $k$-queues $G$ on $n$ vertices such that there is a partition $\{E_1,E_2,\dots,E_k\}$ of $E(G)$, and each $G[E_i]$ is a queue with $|E_i|=m_i$. Then 
\begin{align*}
g(n;m_1,m_2,\dots,m_k)\leq\prod_{i=1}^kg(n,m_i).
\end{align*}
Now $m_1\leq\frac{n}{2}$, as otherwise $m>\frac{kn}{2}\geq m$. Let $j$ be the maximum index such that $m_j\leq\frac{n}{2}$. By \lemref{NumberRainbowEdges},
\begin{align*}
g(n;m_1,m_2,\dots,m_k)
&\leq
\bracket{\prod_{i=1}^j\binom{n}{2m_i}c^{2m_i}}
\bracket{c^n}^{k-j}.
\end{align*}
Now $\sum_{i=1}^jm_i\leq m-\half(k-j)n$. Thus
\begin{align*}
g(n;m_1,m_2,\dots,m_k)
&\leq
\bracket{\prod_{i=1}^j\binom{n}{2m_i}}
\bracket{c^{2m-(k-j)n}}
\bracket{c^{(k-j)n}}\\
&\leq
c^{2m} \prod_{i=1}^k\binom{n}{2m_i}.
\end{align*}
We can suppose that $k$ divides $2m$. It follows (see \citep{BMW-EJC06}) that 
\begin{align*}
g(n;m_1,m_2,\dots,m_k)
\leq
c^{2m}\binom{n}{2m/k}^k.
\end{align*}
It is well known \citep[Proposition~1.3]{Jukna01} that $\binom{n}{t}<(\e n/t)^t$. Thus with $t=2m/k$ we have
\begin{align*}
g(n;m_1,m_2,\dots,m_k)
&<
\bfrac{c\e kn}{2m}^{2m}.
\end{align*}
Clearly
\begin{align*}
g(n,m,k)\leq\sum_{m_1,\dots,m_k}g(n;m_1,m_2,\dots,m_k),
\end{align*}
where the sum is taken over all nonnegative integers $m_1\leq m_2\leq\dots\leq m_k$ such that $\sum_im_i=m$. The number of such partitions \citep[Proposition~1.4]{Jukna01} is at most 
$$\binom{k+m-1}{m}<\binom{2m}{m}<2^{2m}.$$
Hence
\begin{align*}
g(n,m,k)\leq 2^{2m}\bfrac{c\e\,kn}{2m}^{2m}.
\end{align*}
\end{proof}

Every ordered graph on $n$ vertices is isomorphic to at most $n!$ labelled graphs on $n$ vertices. Thus \lemref{NumberRainbowEdges} has the following corollary.

\begin{corollary}
\corlabel{NumberQueue}
For all integers $k$ such that $\frac{2m}{n}\leq k\leq m$, the number of labelled $n$-vertex $m$-edge graphs with queue-number at most $k$ is at most $$\bfrac{ckn}{m}^{2m}n!,$$ for some absolute constant $c$.\qed
\end{corollary}


\begin{proof}[Proof of \thmref{QueueDegree}] 
Let $k$ be the minimum integer such that every simple $\Delta$-regular graph with $n$ vertices has queue-number at most $k$. Thus the number of labelled simple $\Delta$-regular graphs on $n$ vertices is at most the number of labelled graphs with $n$ vertices, $\half\Delta n$ edges, and queue-number at most $k$. By \lemref{NumberRegular} and \corref{NumberQueue},
$$\bfrac{n}{3\Delta}^{\Delta n/2}
\leq 
\bfrac{ck}{\Delta}^{\Delta n}n!
\leq 
\bfrac{ck}{\Delta}^{\Delta n}n^n.$$
Hence 
$k\geq\sqrt{\Delta}n^{1/2-1/\Delta}/(\sqrt{3}c)$.
\end{proof}


\def\soft#1{\leavevmode\setbox0=\hbox{h}\dimen7=\ht0\advance \dimen7
  by-1ex\relax\if t#1\relax\rlap{\raise.6\dimen7
  \hbox{\kern.3ex\char'47}}#1\relax\else\if T#1\relax
  \rlap{\raise.5\dimen7\hbox{\kern1.3ex\char'47}}#1\relax \else\if
  d#1\relax\rlap{\raise.5\dimen7\hbox{\kern.9ex \char'47}}#1\relax\else\if
  D#1\relax\rlap{\raise.5\dimen7 \hbox{\kern1.4ex\char'47}}#1\relax\else\if
  l#1\relax \rlap{\raise.5\dimen7\hbox{\kern.4ex\char'47}}#1\relax \else\if
  L#1\relax\rlap{\raise.5\dimen7\hbox{\kern.7ex
  \char'47}}#1\relax\else\message{accent \string\soft \space #1 not
  defined!}#1\relax\fi\fi\fi\fi\fi\fi} \def\cprime{$'$}

\end{document}